\documentclass[a4paper,10pt]{article}
\usepackage{amssymb}
\usepackage[babel=true,kerning=true]{microtype}

\usepackage{amsthm,leqno}

\usepackage{amsfonts,amsmath,amssymb,enumerate}
\usepackage{graphicx}
\usepackage{tikz}
\usetikzlibrary{decorations.markings}
\usetikzlibrary{plotmarks}
\usetikzlibrary{patterns}

\newcommand{\ds}{\displaystyle}

\newcommand{\be}{\begin{equation*}}
\newcommand{\ee}{\end{equation*}}
\newcommand{\beq}{\begin{equation}}
\newcommand{\eeq}{\end{equation}}
\newcommand{\begincal}{\begin{eqnarray*}}
\newcommand{\fincal}{\end{eqnarray*}}

\newtheorem{thm}{Theorem}[section]
\newtheorem{lemma}{Lemma}[section]
\newtheorem{cor}{Corollary}[section]
\newtheorem{prop}{Proposition}[section]
\newtheorem{defi}{Definition}[section]

\newcommand{\eps}{{\varepsilon}}
\newcommand{\R}{{\mathbb R}}

\newcommand{\N}{{\mathbb N}}

\newcommand{\Z}{{\mathbb Z}}

\def\eps{\varepsilon}
\def\ds{\displaystyle}

\title{Energy quantization for biharmonic maps}
\author{Paul Laurain \& Tristan Rivi\`ere}
\date{december 2011}

\begin{document}
\maketitle
\begin{abstract}
In the present work we establish an energy quantization (or energy identity) result for solutions to scaling invariant variational problems in dimension $4$ which includes   biharmonic maps (extrinsic and intrinsic). To that aim we first establish an angular energy quantization for solutions to critical linear 4th order elliptic systems
with antisymmetric potentials. The method is inspired by the one introduced by the authors previously in \cite{LaR} for 2nd order problems.
\end{abstract}
\section*{Introduction}
Let $N$ be a $C^3$-submanifold of $\R^k$. Let $B_1$ the unit ball of $\R^n$ and  $u\in W^{1,2}(B_1,N)$ then we can define the Dirichlet energy of $u$ as
$$D(u)=\frac{1}{2}\int_{B_1} \vert \nabla u\vert^2\, dx .$$
The critical point of $D$ are the so called harmonic maps for which there have been an extensive theory developed. In particular, when $n=2$ since in that case the functional is conformally invariant, it has been proved that the harmonic maps have some special properties, in particular an energy quantization for sequences of bounded energy, see \cite{Pa} for instance.\\

In this paper, we consider still quadratic scaling invariant problems but in dimension $n=4$ this time. In that case, there is several way to define an equivalent of the Dirichlet functional. Since we look for a scaling invariant quadratic functional the gradient has to be replaced by some expression involving second derivatives. The simplest example 
is given by
$$ E(u)= \frac{1}{4} \int_{B_1} \vert \Delta u\vert^2\, dx .$$
The critical point of this functional are called {\it extrinsic biharmonic maps}. The term extrinsic comes from the fact that this functional (and consequently its critical points) depends on the choice of the embedding of $N$ into $\R^k$. Trying to remedy to this lack of intrinsic nature of the problem, one can instead  consider the following functional
$$ I(u)= \frac{1}{4} \int_{B_1} \vert (\Delta u)^T\vert^2\, dx ,$$
where $(\Delta u)^T$ is the projection of $\Delta u$ onto $T_u N$. The critical point of $I$ will be called {\it intrinsic biharmonic maps}. We can further introduce other functionals
of similar nature and we refer to  \cite{Mo} for more examples. The Euler Lagrange equation satisfied  by the biharmonic maps have been computed in particular in \cite{Wa2}.
One shows that $u\in W^{2,2}(B_1,N)$ is an extrinsic (resp.  intrinsic)  biharmonic map if and only if  $u$ satisfies
$$ T_e(u)\equiv \Delta^2 u - \Delta  (B(u)(\nabla u, \nabla u)) - 2 \nabla\cdot \langle \Delta u, \nabla P(u)\rangle +  \langle \Delta (P(u)),\Delta u\rangle=0,$$
respectively 
\be
\begin{split}
 T_i(u)&\equiv \Delta^2 u - \Delta  (B(u)(\nabla u, \nabla u)) - 2 \nabla\cdot \langle \Delta u, \nabla P(u)\rangle + \langle \Delta (P(u)),\Delta u\rangle \\
 &-P(u)\left(B(u) (\nabla u, \nabla u )  \nabla_u B(u) (\nabla u, \nabla u ) \right)- 2 B(u)(\nabla u , \nabla u) B(u) (\nabla u, \nabla P(u))=0,
 \end{split}
 \ee
where $P$ and $B$ are the orthogonal projection onto $T_u N$ and the second fundamental form of $N$\footnote{see section \ref{s1} for precise definitions}. Since our result applies indistinctly to extrinsic as well as to intrinsic biharmonic maps, except when it is necessary,  in what follow {\bf we will indifferently employ the denomination biharmonic map  for both extrinsic biharmonic map and intrinsic biharmonic map}. We observe that these equations are of the form,
$$\Delta^2  u = \sum_{\substack{ \alpha_1+\cdots+\alpha_4=4\\ 0\leq \alpha_i<4}} c_\alpha(u) \,  \partial^{\alpha_1}  u\,\partial^{\alpha_2}  u\, \partial^{\alpha_3}  u\, \partial^{\alpha_4}  u ,$$
which make them critical in the sense that there is no hope to use directly standard $L^p$-theory for proving regularity or compactness results. The critical nature of an elliptic problem is
characterized by possible loss of compactness at isolated points. In order to fully describe this concentration-compactness phenomenon
one has to understand ''how much'' energy is lost at these isolated points. {\it Energy quantization} means that the energy lost corresponds exactly to the sum of 
the energies of the so called {\it bubbles} - or rescaled elementary solutions on $S^4$ - concentrating at these points. The word {\it quantization} refers to the fact that
the bubbles cannot have arbitrary small energy and in some problems it is even known that they can realize only a discrete set of values.

Our main result in this paper is the energy quantization result for biharmonic maps. In fact we are proving something stronger considering more generally sequences of approximate solutions of biharmonic maps. To that aim  we need the following definition. 

\begin{defi} Let $N$ be a $C^3$-submanifold of $\R^k$, $p\geq1$, $f \in L^p(B_1,\R^k)$  and $u\in W^{2,2}(B_1,N)$. $u$ is $f$-approximate biharmonic maps if $u$ satisfies
$$ T_i (u)=f \hbox{ or } T_e (u)=f .$$ 
\end{defi}

Hence, we are in position to sate our main result.

\begin{thm}
\label{th1}
Let $N$ be a $C^3$-submanifold of $\R^k$, $p>1$, $f_n \in L^p(B_1,\R^k)$ and $u_n\in W^{2,2}(B_1,N)$ be a sequence of  $f_n$-approximate biharmonic maps with bounded energy, i.e.
\beq
\label{be}
\int_{B_1} \left(\vert \nabla^2 u_n \vert^2 + \vert \nabla u_n \vert^4  + \vert f_n\vert^p\right)\, dz  \leq M.
\eeq
Then there exists $f\in L^{p}(B_1,\R^k)$,  $u_\infty\in W^{2, 1}(B_1,N)$ a $f$-approximate  bihamonic map and
\begin{enumerate}[i)]
\item $\omega^1,\dots,\omega^l$ some biharmonic maps of $\R^4$ to $N$,
\item $a_n^1,\dots,a_n^l$ a family of converging sequences of points of $B_1$,
\item  $\lambda_n^1,\dots,\lambda_n^l$ a family of  sequences of positive reals converging all to zero,
\end{enumerate}
such that, up to a subsequence,  
$$u_n\rightarrow u_\infty \hbox{ in } W^{2,q}_{loc}(B_1 \setminus \{ a^1_\infty,\dots,a^l_\infty\}) \hbox{ for all } q < \frac{2p}{2-p}$$
and 
$$ \left\Vert \nabla^2  \left( u_n -u_\infty -\sum_{i=1}^l \omega^i_n \right) \right\Vert_{L^2_{loc}(B_1)}+ \left\Vert \nabla  \left( u_n -u_\infty -\sum_{i=1}^l \omega^i_n \right) \right\Vert_{L^4_{loc}(B_1)} \rightarrow 0,$$
where $\omega^i_n=\omega^i(a_n^i+\lambda_n^i \, . \, )$. Moreover, if $N$ is $C^{l+3}$ and $f_n$ is bounded in $C^{l,\eta}(B_1,\R^k)$ then the convergence of $u_n$ to $u_\infty$ is in  $C^{l+4,\nu}(B_1 \setminus \{ a^1_\infty,\dots,a^l_\infty\})$ for any $0\leq \nu <\eta$.
\end{thm}
Observe that sequences of biharmonic maps into a smooth manifold holds in $C^{\infty}_{loc}$. Such a result was already know for intrinsic biharmonic maps, see \cite{HP1} and \cite{HP2}, or for extrinsic biharmonic maps into a sphere, see \cite{Wa4}. Here, the method employed seems particularly robust since, not only it is identical for both extrinsic and intrinsic biharmonic maps but it applies moreover to a larger class of scaling invariant problems. As an illustration of this fact we prove that the method applies to 
the following general lagrangians
\beq
\label{gl}
\int_{B_1} \left(\vert \Delta u \vert^2\, dx  +u^* \Omega\right)  \hbox{ or } \int_{B_1}\left( \vert (\Delta u)^T \vert^2 \, dx  +u^* \Omega\right) ,
\eeq
where $\Omega$ is an arbitrary smooth 4-form of $\R^k$.

The method we use goes first through the proof of an angular energy quantization result\footnote{see the end of  section \ref{s4} for a precise statement} for sequences  of solutions to the general critical 4th order elliptic system system with antisymmetric potentials
introduced by Lamm and Riviere \cite{LR} . We follow in fact the approach that we originally introduced in \cite{LaR} for second order problems.
We have good reasons to think that the method could further be extended for proving a general energy quantization result for polyharmonic maps in critical dimension, (see the $\eps$-regularity for polyharmonic maps in \cite{GSZG} and \cite{GasSch}  for the general case, see also \cite{Ru}).\\

As an immediate consequence  of theorem \ref{th1}, we get the asymptotic behavior of biharmonic maps flow. A weak solution to the extrinsic biharmonic map flow is a map $u\in W^{2,2}([0,+\infty [\times B_1,N)$ satisfying
\beq
\label{bf}
\left\{
\begin{split}
& \frac{\partial u}{\partial t} + \Delta^2 u = \Delta  (B(u)(\nabla u, \nabla u)) + 2 \nabla\cdot \langle \Delta u, \nabla P(u)\rangle - \langle \langle \Delta (P(u)),\Delta u\rangle \hbox { on } [0,+\infty[ \times B_1 \\
&u= u_0 \hbox { on } \{0\} \times B_1
\end{split}
\right.
\eeq
where $u_0 \in W^{2,2}(B_1,N)$.  Several existence results have been established for (\ref{bf}), see for instance \cite{Lamm}  for small initial data or \cite{Gas} and \cite{Wa3} for solution with finitely many singular times and arbitrary initial data.  All these solutions satisfy the following energy identity
\beq
\label{ei} 
2\int_{0}^T \int_{B_1} \left\vert \frac{\partial u}{\partial t} \right\vert^2 \, dxdt + \int_{B_1} \vert \Delta u\vert^2 \, dx \leq   \int_{B_1} \vert\Delta u_0\vert^2 \, dx  \hbox{ for all } T\geq 0.
\eeq

\begin{cor}
\label{th3}
Let $N$ be a $C^3$-submanifold of $\R^k$ and $u_0\in W^{2,1}(B_1,N)$ and  $u\in W^{2,2}([0,+\infty [\times B_1,N)$ be a global solution of (\ref{bf}) satisfying the energy inequality (\ref{ei}). Then there exist $t_n $ a sequence of positive real such that $t_n\rightarrow +\infty$, a biharmonic map  $u_\infty\in W^{2, 1}(B_1,N)$, $l\in\N$,  $\omega^1,\dots,\omega^l$ some biharmonic maps of $\R^4$ to $N$ and  $a_n^1,\dots,a^l_n$ a family of  of points of $B_1$ converging to
$a_\infty^1,\dots,a^l_\infty$ , such that 
$$u(t_n, \, . \,) \longrightarrow u_\infty \hbox{ on } W^{2,p}_{loc}(B_1 \setminus \{ a^1_\infty,\dots,a^l_\infty\}) \hbox{ for all } p\geq 1 $$
and 
$$ \left\Vert \nabla^2  \left( u(t_n, \, . \,)  -u_\infty -\sum_{i=1}^l \omega^i_n \right) \right\Vert_{L^2_{loc}(B_1)}+ \left\Vert \nabla  \left( u(t_n, \, . \,)  -u_\infty -\sum_{i=1}^l \omega^i_n \right) \right\Vert_{L^4_{loc}(B_1)} \rightarrow 0,$$
where $\omega^i_n=\omega^i(a_n^i+\lambda_n^i \, . \, )$.
\end{cor}
In fact, thanks to (\ref{ei}), we easily prove that there exists $t_n$ such that $u(t_n, \, .\,)$ satisfies the hypothesis of theorem \ref{th1} with $p=2$.\\

The paper is organized as follows: in section \ref{s1}, we rewrite the equations in order to apply the theory of Lamm and Riviere, in section \ref{s2} we recall the main results of Lamm and Riviere and we prove an $\eps$-regularity result for biharmonic maps, in section \ref{s3} we derive the key estimate in Lorentz space for the angular derivatives  in a annular region of arbitrary conformal type, finally in section \ref{s4} we prove our main result postponing technical result to section \ref{s5} and \ref{sp}.

\medskip {\bf Acknowledgements} : This work was initiated as the first author was visiting  the {\it Forschungsinstituts f\"ur Mathematik} at E.T.H. (Zurich). He would like to thank 
the institute for its hospitality and the excellent working conditions.

\section{Biharmonic equation in normal form}
\label{s1}

Let $N\subset \R^k$ be a $C^3$-submanifold, there exists  $\delta>0$ such that $\Pi : N_\delta \rightarrow N$, the nearest projection map, is well defined and $C^3$, where $N_\delta=\{y\in\R^k\,\vert\, d(y,N)\leq \delta\}$. Let, for $y\in N$,  $ P(y)\equiv \nabla \Pi (y):\R^k \rightarrow T_y N$ be the orthogonal projection, and $ P^\bot(y)\equiv Id-\nabla\Pi (y):\R^k \rightarrow (T_y N)^\bot$. In the following, we will write $P$(resp. $P^\bot$) instead of  $P(y)$(resp. $P^\bot(y)$) and we will identify these linear transformations with their matrix representations in $\mathcal{M}_k$. We also note that these projections are in $W^{2,2}(B_1,\mathcal{M}_k)$ as soon as $u$ is in  $W^{2,2}(B_1,N)$ . Finally, let $B(\, .\,)(\, . \, , \, . \,)$ be the second fundamental form of $N\subset\R^k$, which is defined by
$$B(y)(Y,Z)=D_Y P^\bot (y) (Z), \; \forall y\in N,Y,Z \in T_y N.$$
We know that, see \cite{Wa1}, that $u\in W^{2,2}(B_1,N)$ is an extrinsic biharmonic map if and only if 
$$ \Delta^2 u\,  \bot\, T_u N \hbox{ almost everywhere},$$
which can be rewritten as follows
\beq
\label{m1} 
\begin{split}
\Delta^2 u &= P^\bot  \Delta^2 u \\
&=div(P^\bot  \nabla \Delta u)-\nabla P^\bot \nabla \Delta u.
\end{split}
\eeq
Then we rewrite the second term of the right hand side as follows
\beq
\label{m2} 
\begin{split}
\nabla P^\bot \nabla \Delta u&= \nabla P^\bot P^\bot \nabla \Delta u + \nabla P^\bot P \nabla \Delta u\\
&= \nabla P^\bot P^\bot  \nabla \Delta u -  P^\bot \nabla P \nabla \Delta u\\
&= 2 \nabla P^\bot P^\bot  \nabla \Delta u+  (\nabla P P^\bot -  P^\bot \nabla P) \nabla \Delta u .
\end{split}
\eeq
But 
\beq
\label{m3} 
\begin{split}
2 \nabla P^\bot P^\bot  \nabla \Delta u& = 2 \nabla P^\bot P^\bot  \nabla \Delta u - 2 \nabla P^\bot \nabla div (P^\bot  \nabla u)\\
&= - 2 \nabla P^\bot \nabla  P^\bot   \Delta u + 2 div( \nabla P^\bot (\nabla P^\bot  \nabla u)) - 2 \Delta P^\bot \nabla P^\bot \nabla u
\end{split}
\eeq
Thanks to (\ref{m1}), (\ref{m2}) and (\ref{m3}), we get  
\be
\begin{split}
\Delta^2 u &= div ( P^\bot  \nabla \Delta u)-div (2\nabla P^\bot (\nabla P^\bot  \nabla u)) \\
&+ 2 \nabla P^\bot \nabla P^\bot  \Delta u + 2\Delta P^\bot \nabla P^\bot \nabla u\\
&-(\nabla P  P^\bot -P^\bot \nabla P) \nabla \Delta u\\
&= \Delta( P^\bot  \Delta u)-div (\nabla P^\bot \Delta u +2\nabla P^\bot (\nabla P^\bot  \nabla u)) \\
&+ 2 \nabla P^\bot \nabla P^\bot  \Delta u + 2\Delta P^\bot \nabla P^\bot \nabla u\\
&-(\nabla P  P^\bot -P^\bot \nabla P) \nabla \Delta u,
\end{split}
\ee
which finally gives the equation of extrinsic biharmonic maps
\beq
\label{ebe}
\begin{split}
\Delta^2 u&= -\Delta (\nabla P^\bot  \nabla u)-div (\nabla P^\bot  \Delta u) \\
&+ 2 \nabla P^\bot \nabla(\nabla P^\bot  \nabla u) + 2 \nabla P^\bot \nabla P^\bot \Delta u\\
&-(\nabla P  P^\bot -P^\bot \nabla P) \nabla \Delta u.
\end{split}
\eeq
For intrinsic biharmonic maps, we need to add some tangent terms, see \cite{Wa2} for details, which gives
\beq
\label{ibe}
\begin{split}
\Delta^2 u&= -\Delta (\nabla P^\bot  \nabla u)-div (\nabla P^\bot  \Delta u) \\
&+ 2 \nabla P^\bot \nabla(\nabla P^\bot  \nabla u) + 2 \nabla P^\bot \nabla P^\bot \Delta u\\
&-(\nabla P  P^\bot -P^\bot \nabla P) \nabla \Delta u\\
&+P\left( \nabla P^\bot  \nabla u \nabla(\nabla P^\bot  \nabla u)\right)\\
&+ 2 \nabla P^\bot  \nabla u \nabla P^\bot  \nabla P .
\end{split}
\eeq

\begin{prop}
\label{pn}
The  equation (\ref{ebe}) and  (\ref{ibe}) can be rewritten in the form
\beq
\label{nf}
\Delta^2 u = \Delta (V \nabla u) + div( w \nabla u) + \nabla \omega \nabla u + F\nabla u ,
\eeq
where $V \in W^{1,2}(B_1, \mathcal{M}_k \otimes \Lambda^{1} \R^4)$, $w\in L^2(B_1, \mathcal{M}_k)$, $\omega \in L^2(B_1, \mathrm{s}o_k)$ and $F \in L^2\cdot W^{1,2}(B_1,  \mathcal{M}_k \otimes \Lambda^{1} \R^4)$ with
\beq
\label{cc}
\begin{split}
&\vert V\vert  \leq C\left(  \vert \nabla u\vert \right)\\
& \vert F \vert  \leq C\left( \left(\vert \nabla^2 u\vert + \vert \nabla u\vert^2 \right) \vert\nabla u \vert \right) \hbox{ almost everywhere },\\
&\vert w \vert + \vert \omega \vert  \leq C\left( \vert \nabla^2 u\vert + \vert \nabla u\vert^2 \right ) 
\end{split}
\eeq
where $C$ is a positive constant which depends only on $N$ .
\end{prop}

{\it Proof of proposition~\ref{pn}:}\\

We give a proof for equation (\ref{ebe}), the intrinsic case will follow easily.\\

From the one hand, we proceed to the following Hodge decomposition
$$  dP  P^\bot -P^\bot d P = d \alpha + d^* \beta,$$
where $ \alpha \in W^{1,2}(B_1, \mathrm{s}o_k)$, $\beta \in W^{1,2}_0(B_1, \Lambda^2 (R^4) \otimes\mathcal{M}_k)$. Hence $\alpha$ and $\beta$ satisfy
$$\Delta \alpha= \Delta P  P^\bot -P^\bot \Delta P,$$
and 
$$\Delta\beta= d P\wedge d P^\bot - d P^\bot \wedge d P .$$
Then $\alpha \in W^{2,2}(B_1,  \mathrm{s}o_k)$, $d^* \beta \in W^{2,(\frac{4}{3},1)}_0(B_1, \Lambda^2 (R^4) \otimes\mathcal{M}_k)$ and we get
\be
\begin{split}
(\nabla P  P^\bot -P^\bot \nabla P) \nabla \Delta u &= d \Delta \alpha \nabla u + \Delta  d^* \beta \nabla u + \Delta((\nabla P  P^\bot -P^\bot \nabla P)\nabla u)\\
& -2 div (\nabla(\nabla P  P^\bot -P^\bot \nabla P)\nabla u)\\
&=\nabla \omega_1 \nabla u + F_1 \nabla u + \Delta (V_1 \nabla u) + div (w_1 \nabla u),
\end{split}
\ee
with $\omega_1 \in L^2(B_1,\mathrm{s}o_k)$, $F_1\in L^{2}\cdot W^{1,2}(B_1,  \mathcal{M}_k \otimes \Lambda^{1} \R^4)$, $V_1\in W^{1,2}(B_1, \mathcal{M}_k \otimes \Lambda^{1} \R^4)$ and $w_1\in L^2(B_1, \mathcal{M}_k)$.\\

From the other hand, we have
$$  2 \nabla P^\bot \nabla(\nabla P^\bot  \nabla u) = F_2 \nabla u,$$
with $F_2 ^l=2 \frac{\partial P^\bot}{\partial y^l} \nabla(\nabla P^\bot  \nabla u) \in L^2\cdot W^{1,2}(B_1,  \mathcal{M}_k \otimes \Lambda^{1} \R^4)$ and
$$2 \nabla P^\bot \nabla P^\bot \Delta u = F_3 \nabla u,$$
with $F_3 ^l=2 \frac{\partial P^\bot}{\partial y^l} \nabla P^\bot  \Delta u \in L^2\cdot W^{1,2}(B_1,  \mathcal{M}_k \otimes \Lambda^{1} \R^4)$, which achieves the proof.\hfill$\square$\\

For general Lagrangian of the form (\ref{gl}), the equation becomes,

$$ T_e(u) = H\left(\frac{\partial u}{\partial x_1},\frac{\partial u}{\partial x_2},\frac{\partial u}{\partial x_3},\frac{\partial u}{\partial x_4}\right) \hbox{ or } T_e(u) =H\left(\frac{\partial u}{\partial x_1},\frac{\partial u}{\partial x_2},\frac{\partial u}{\partial x_3},\frac{\partial u}{\partial x_4}\right),$$
where $H$ is the $4$-form on $\R^k$ into $\R^k$ defined by
$$d\Omega(U,V,W,X,Y)= U_iH^i(V,W,X,Y) \hbox{  for all  } U,V,W,X,Y \in \R^k .$$
Hence we have 
$$H\left(\frac{\partial u}{\partial x_1},\frac{\partial u}{\partial x_2},\frac{\partial u}{\partial x_3},\frac{\partial u}{\partial x_4}\right) =F\nabla u,$$
with $F \in L^2\cdot W^{1,2}(B_1, \mathcal{M}_k \otimes\Lambda^1\R^4 )$.

\section{Preliminaries}
\label{s2}
First, we recall the main result of \cite{LR} that provides a divergence form  to elliptic 4th order system of the kind (\ref{nf}) under small energy assumption. This will be one of the main tools in order to obtain the estimate needed for the energy quantization.

\begin{thm}[Theorem 1.4 \cite{LR}]
\label{TLR0} There exists $\eps>0$ and $C>0$ depending only on $N$, such that the following holds:  Let  $V \in W^{1,2}(B_1, \mathcal{M}_k \otimes \Lambda^{1} \R^4)$, $w\in L^2(B_1, \mathcal{M}_k)$, $\omega \in L^2(B_1,\mathrm{s}o_k)$ and  $F \in L^2\cdot W^{1,2}(B_1,  \mathcal{M}_k \otimes \Lambda^{1} \R^4)$ such that
$$\Vert V \Vert_{W^{1,2}} +\Vert w \Vert_2 + \Vert \omega \Vert_2 + \Vert F \Vert_{L^2 \cdot W^{1,2}} < \eps,$$
then there exists $A\in L^\infty \cap W^{2,2}(B_1, \mathcal{G}l_k)$ and $B\in  W^{1,\frac{4}{3}}(B_1, \mathcal{M}_k\otimes \Lambda^2 \R^4)$ such that 
$$ \nabla \Delta A + \Delta A V -\nabla A w + A (\nabla \omega +F)= \mathrm{curl} B,$$
and
$$\Vert A\Vert_{W^{2,2}} + d(A, \mathcal{SO}_n) + \Vert B\Vert_{W^{1,\frac{4}{3}}} \leq C \left(\Vert V \Vert_{W^{1,2}} + \Vert w\Vert_2 +\Vert \omega \Vert_2 +    \Vert F \Vert_{L^2 \cdot W^{1,2}}\right).$$
\end{thm}

Thanks to the previous theorem, we are in position to rewrite equations of the form (\ref{nf}) in divergence form.
\begin{thm}[Theorem 1.2 and 1.4  \cite{LR}]
\label{TLR} 
There exists $\eps>0$ and $C>0$ depending only on $N$, such that if  $u \in W^{2,2}(B_1,\R^k)$  satisfies
$$\Delta^2 u = \Delta (V \nabla u) + div( w \nabla u) + \nabla \omega \nabla u + F\nabla u + f ,$$ 
where $V \in W^{1,2}(B_1, \mathcal{M}_k \otimes \Lambda^{1} \R^4)$, $w\in L^2(B_1, \mathcal{M}_k)$, $\omega \in L^2(B_1,\mathrm{s}o_k)$, $F \in L^2\cdot W^{1,2}(B_1,  \mathcal{M}_k \otimes \Lambda^{1} \R^4)$ and  $f \in L^1(B_1,\R^k)$  with
$$\Vert V \Vert_{W^{1,2}} +\Vert w \Vert_2 + \Vert \omega \Vert_2 + \Vert F \Vert_{L^2 \cdot W^{1,2}} < \eps,$$
then there exists $A\in L^\infty \cap W^{2,2}(B_1, \mathcal{G}l_k)$ and $B\in  W^{1,\frac{4}{3}}(B_1, \mathcal{M}_k\otimes \Lambda^2 \R^4)$ such that 
$$\Vert A\Vert_{W^{2,2}} + d(A, \mathcal{SO}_n) + \Vert B\Vert_{W^{1,\frac{4}{3}}} \leq C \left(\Vert V \Vert_{W^{1,2}} + \Vert w\Vert_2 +\Vert \omega \Vert_2 +    \Vert F \Vert_{L^2 \cdot W^{1,2}}\right)$$
and
$$\Delta(A \Delta u)= div\left(2\nabla A \Delta u-\Delta A \nabla u + A w \nabla u + \nabla A(V\nabla u)-A \nabla (V \nabla u) -B\nabla u\right) + Af.$$
\end{thm}

A first consequence of  the previous theorem, is the $\eps$-regularity for biharmonic maps. It can also be compared with the corresponding result
established for second order problems in theorem 3.2 of \cite{LaR}.  

\begin{thm}
\label{ereg}
Let $p>1$, there exists $\eps>0$ and $C_{p}>0$ such that if $u\in W^{2,2}(B_1, \R^k)$ (resp. $u\in W^{2,2}(\R^4, \R^k)$), $f\in L^p (B_1, \R^k)$ (resp. $f\in L^p (\R^4, \R^k)$), $V \in W^{1,2}(B_1, \mathcal{M}_k \otimes \Lambda^{1} \R^4)$(resp.  $V \in W^{1,2}(\R^4, \mathcal{M}_k \otimes \Lambda^{1} \R^4)$), $w\in L^2(B_1, \mathcal{M}_k)$ (resp. $w\in L^2(\R^4, \mathcal{M}_k)$), $\omega \in L^2(B_1, \mathrm{s}o_k)$(resp.  $\omega \in L^2(\R^4, \mathrm{s}o_k)$) and $F \in L^2\cdot W^{1,2}(B_1,  \mathcal{M}_k \otimes \Lambda^{1} \R^4)$(resp. $F \in L^2\cdot W^{1,2}(\R^4,  \mathcal{M}_k \otimes \Lambda^{1} \R^4)$) satisfy (\ref{cc})  and 
$$\Vert \nabla^2 u\Vert_2+ \Vert \nabla u\Vert_4  \leq \eps,$$
then
\begin{enumerate}
\item($\eps$-regularity) If  $u\in W^{2,2}(B_1, \R^k)$ is a solution of 
$$\Delta^2 u = \Delta (V \nabla u) + div( w \nabla u) + \nabla \omega \nabla u + F\nabla u + f \hbox{ on }  B_1,$$ 
 then we have
$u\in W^{2, \bar{p} }(B_\frac{1}{2}, \R^k)$,  where $\bar{p} =\frac{2p}{p-2}$ if $p<2$ else any $\bar{p}\geq 2$  and 
$$ \Vert \nabla^2 u\Vert_{L^{\bar{p}}\left(B_{\frac{1}{2}}\right)}  +  \Vert  \nabla u\Vert_{L^{2\bar{p}}\left(B_{\frac{1}{2}}\right)}  \leq C_{p}   \left( \Vert \nabla^2 u\Vert_{L^2(B_1)}+ \Vert \nabla u\Vert_{L^4(B_1)}  + \Vert f\Vert_p \right).$$
Moreover, if $N$ is smooth and $f\in C^{l,\eta}$ for $l\in\N$ and $\eta>0$ then we can replace $W^{4, \bar{p}}$ by $C^{l+4,\eta}$. 
\item(Energy gap) If  $u\in W^{2,2}(\R^4, \R^k)$ is a solution of 
$$\Delta^2 u = \Delta (V \nabla u) + div( w \nabla u) + \nabla \omega \nabla u + F\nabla u \hbox{ on } \R^4,$$
then $u$ is identically equal to zero. 
\end{enumerate} 
\end{thm}

The proof  of theorem \ref{ereg} could be achieved almost following lemma 3.1 of \cite{LR}. We give however an independent proof of this fact that shed new lights on the problem.\\

{\it Proof  of theorem \ref{ereg}:}\\

Let $0<\eps<1$ such that, thanks to (\ref{cc}), hypothesis of  theorem \ref{TLR} are satisfied. Then we can rewrite  our equation as
$$ \Delta (A\Delta u)= div(K) + Af,$$
where $A\in L^\infty \cap W^{2,2}(B_1, \mathcal{G}l_k)$ and $K\in L^2\cdot W^{1,2} \subset L^{\frac{4}{3},1}$ satisfy
$$\Vert A\Vert_{W^{2,2}} + d(A, \mathcal{SO}_n) + \Vert K\Vert_{L^{\frac{4}{3},1}} \leq C \left(\Vert \nabla^2 u\Vert_2+ \Vert \nabla u\Vert_4+ \Vert V \Vert_{W^{1,2}} +\Vert w \Vert_2 + \Vert \omega \Vert_2 + \Vert F \Vert_{L^2 \cdot W^{1,2}}\right)$$
where $C$ is independent of $u$.\\

 Let $p\in B_{\frac{1}{2}}$ and $0<\rho<\frac{1}{2}$. We decompose $A\Delta u$ on $B_\rho(p)$ as  $A\Delta u= C+ D$ where $C\in W_0^{1,2}(B_\rho(p))$ and $D\in W^{1,2}(B_\rho(p))$. Then  $C$  satisfies
$$\Delta C = div(K) + Af \hbox{ on } B_\rho(p)$$
and $D$ satisfies
$$
\Delta D = 0 \hbox{ on } B_\rho(p).
$$
Thanks to the standard $L^p$-theory and Sobolev embeddings, we get
\beq
\label{EC}
\left(\int_{B_\rho(p)} \vert C\vert^2 \, dx\right)^\frac{1}{2}  \leq C \left(\Vert K\Vert_{\frac{4}{3}}+ \rho^\frac{4(p-1)}{p} \Vert f\Vert_p\right)  \leq C \left(\eps  \Vert \nabla^2 u \Vert_{2} + \frac{\eps}{\rho} \Vert \nabla u \Vert_2 +  \rho^\frac{4(p-1)}{p} \Vert f\Vert_p\right),
\eeq
where $C$ is a positive constant in dependent of $u$.\\

Using the fact that $D$ is harmonic, we have that $\delta \mapsto \frac{1}{(\delta\rho)^4} \int_{B_{\delta\rho}(p)} \vert D\vert^2\, dx$  is an increasing function and hence for all $\delta\in]0,1[$
we deduce,
\beq
\label{EC2}
\int_{B_{\delta \rho}(p)}  \vert D\vert^2\, dx  \leq \delta^4 \int_{B_\rho(p)}  \vert D\vert^2\, dx. 
\eeq
We then decompose $u$ as follows :  $u= E+ F$ where $E\in W_0^{1,4}(B_\rho(p))$ and $F\in W^{1,4}(B_\rho(p))$ satisfy
$$
\Delta E = A^{-1}(C+D)\hbox{ on } B_\rho(p)$$
and $F$ satisfies
$$
\Delta F = 0 \hbox{ on } B_\rho(p) .$$
Thanks to the standard $L^p$-theory and Sobolev embeddings, we get
\beq
\label{EC3}
\frac{1}{\rho} \left(\int_{B_\rho(p)} \vert \nabla E\vert^2 \, dx\right)^\frac{1}{2}  \leq C \left( \left(\int_{B_\rho(p)}\vert C\vert^2\, dx\right)^\frac{1}{2} +\left(\int_{B_\rho(p)}\vert D\vert^2\, dx  \right)^\frac{1}{2}\right).
\eeq
where $C$ is a positive constant in dependent of $u$.\\

The function  $\delta \mapsto \frac{1}{(\delta\rho)^4} \int_{B_{\delta\rho}(p)} \vert \nabla F \vert^2\, dx$  is increasing since $F$ is harmonic and  we have have again, for all $\delta\in]0,1[$,
\beq
\label{EC4}
\frac{1}{(\delta\rho)^2}\int_{B_{\delta\rho}(p)}  \vert \nabla F \vert^2\, dx  \leq \frac{\delta^2}{\rho^2} \int_{B_\rho(p)}  \vert \nabla F \vert^2\, dx.
\eeq
Then, thanks to (\ref{EC}), (\ref{EC2}), (\ref{EC3}) and (\ref{EC4}), for  $\delta$ and $\eps$ small enough (with respect to some constant independent of $u$), we have 
$$\int_{B_{\delta \rho}(p)} \left(\vert \nabla^2 u \vert^2+\frac{1}{(\delta\rho)^2}\vert \nabla u \vert^2\right) \, dx \leq \frac{1}{2} \int_{B_{\rho}(p)} \left(\vert \nabla^2 u \vert^2 +\frac{1}{\rho^2} \vert \nabla u \vert^2\right) \, dx + C\delta^\frac{4(p-1)}{p}  \Vert f\Vert_p^2.$$

Iterating this inequality gives the following Morrey type estimate : there exists $\alpha>0$ and $C>0$ such that 
$$\sup_{p\in B_\frac{1}{2} , 0<\rho <\frac{1}{2}} \rho^{-\alpha}\left(\int_{B_{\rho}(p)} \left(\vert \nabla^2 u \vert^2+\frac{1}{\rho^2}\vert \nabla u \vert^2\right) \, dx\right) \leq C \Vert f\Vert_p .$$
Then
$$\sup_{p\in B_\frac{1}{2} , 0<\rho <\frac{1}{2}} \rho^{-\alpha} \int_{B_{\rho}(p)} \vert \Delta^2 u \vert \, dx   \leq C \Vert f\Vert_p. $$
Then a classical estimate on Riesz potentials gives, for all $p\in B_{\frac{1}{3}}$
\be
\begin{split}
&|\Delta u|(p)\le (C \Vert f\Vert_p)\frac{1}{|x|^2}\ast \chi_{B_{\frac{1}{2}}}\ |\Delta^2 u|+C \Vert \nabla^2 u\Vert_{L^2(B_1)},\\
& |\nabla u|(p)\le (C \Vert f\Vert_p)\frac{1}{|x|}\ast \chi_{B_{\frac{1}{2}}}\ |\Delta^2 u|+C \Vert \nabla u\Vert_{L^2(B_1)},
\end{split}
\ee
where $\chi_{B_{\frac{1}{2}}}$ is the characteristic function of the ball $B_{\frac{1}{2}}$. Together with injections proved by Adams in \cite{Ad}, see also 6.1.6 of \cite{Gra1}, the latter shows that
$$
\Vert \nabla^2 u \Vert_{L^{r}\left(B_{\frac{1}{3}}\right)} + \Vert \nabla u \Vert_{L^{2r}\left(B_{\frac{1}{3}}\right)} \leq C  \left( \Vert f\Vert_p  +\Vert \nabla^2 u\Vert_2+ \Vert \nabla u\Vert_4 \right)\quad ,
$$
for some $r>1$. Then bootstrapping this estimate, we get 
$$\Vert  \nabla^2 u\Vert_{L^{\bar{p}} (B_{\frac{1}{4}})}  + \Vert  \nabla u\Vert_{L^{2\bar{p}}(B_{\frac{1}{4}})}  \leq C \left( \Vert f\Vert_p  +\Vert \nabla^2 u\Vert_2+ \Vert \nabla u\Vert_4 \right),$$
where $\bar{p}$ is the limiting exponent of the bootstrapping given by the Sobolev injection of $W^{2,p}$ into $L^{\bar{p}}$ if $p<2$. Indeed, thanks to (\ref{cc}), the only limiting term for the bootstrap is the regularity of $f$.\\

Now, we can easily derive the proof of the energy gap. Indeed, thanks to the previous estimate, we easily see that  for some $q>2$ we get
$$\Vert \nabla^2 u\Vert_{L^q(B_R)} +  \Vert \nabla u\Vert_{L^{2q}(B_R)} \leq C  \frac{\Vert u\Vert_{W^{2,2}}}{R^{2-\frac{4}{q}}} \hbox {for all } R>0,$$
which proves that $u\equiv 0$.\hfill$\square$
\section{Uniform estimate in annular region}
\label{s3}
In this section, we derive a strong estimate for angular derivatives in an annular region independently of the conformal class. 

\begin{thm}
\label{EN} There exist $\eps>0$ and $C>0$ depending only on $k$, such that if $0<r<\frac{1}{4}$, $p>1$ and   $u \in W^{2,2}(B_1\setminus B_r,\R^k)$  satisfies
$$\Delta^2 u = \Delta (V \nabla u) + div( w \nabla u) + \nabla \omega \nabla u + F\nabla u +f ,$$ 
where $V \in W^{1,2}(B_1\setminus B_r, \mathcal{M}_k \otimes \Lambda^{1} \R^4)$, $w\in L^2(B_1\setminus B_r, \mathcal{M}_k)$, $\omega \in L^2(B_1\setminus B_r,\mathrm{s}o_k)$, $F \in L^2\cdot W^{1,2}(B_1\setminus B_r,  \mathcal{M}_k \otimes \Lambda^{1} \R^4)$ and $f\in L^p(B_1,\R^k)$  with
$$\Vert V \Vert_{W^{1,2}} +\Vert w \Vert_2 + \Vert \omega \Vert_2 + \Vert F \Vert_{L^2 \cdot W^{1,2}} < \eps,$$
then 
$$\left\Vert \nabla^T \nabla u \right\Vert_{L^{2,1}\left(B_\frac{1}{4}\setminus B_{4r}\right)} \leq C \left(1 + \Vert \nabla^2 u\Vert_{L^2 (B_1\setminus B_r)} + \Vert \nabla u\Vert_{L^4 (B_1\setminus B_r)} + \Vert f\Vert_{L^p (B_1\setminus B_r)}\right) ,$$
where $\nabla^T f= \nabla f - \frac{\partial f}{\partial r} \frac{\partial }{\partial r}$.  
\end{thm}

{\it Proof of theorem \ref{EN}:}\\

Using some classical extension theorem, we see that there exist $\tilde{V} \in W^{1,2}(B_1, \mathcal{M}_k \otimes \Lambda^{1} \R^4)$, $\tilde{w}\in L^2(B_1, \mathcal{M}_k)$, $\tilde{\omega} \in L^2(B_1,\mathrm{s}o_k)$ and $\tilde{F} \in L^2\cdot W^{1,2}(B_1,  \mathcal{M}_k \otimes \Lambda^{1} \R^4)$ such that $\tilde{V}=V$, $\tilde{w}=w$, $\tilde{\omega}=\omega$ and $\tilde{F}=F$ on $B_1\setminus B_r$ and 
$$\Vert \tilde{V} \Vert_{W^{1,2}} +\Vert \tilde{w} \Vert_2 + \Vert \tilde{\omega} \Vert_2 + \Vert \tilde{F} \Vert_{L^2 \cdot W^{1,2}} <2 \eps,$$

Thanks to theorem \ref{TLR0}, for $0<\eps<\frac{1}{2}$ small enough, there exist $A\in L^\infty \cap W^{2,2}(B_1, \mathcal{G}l_k)$ and  $B \in W^{1,(\frac{4}{3},1)} (B_1)$ such that 
$$ d(A,\mathcal{SO}_k) + \Vert A \Vert_{W^{2,2}}+ \Vert B\Vert_{W^{1,(\frac{4}{3},1)}}  \leq C\left(\Vert \tilde{V} \Vert_{W^{1,2}} +\Vert \tilde{w} \Vert_2 + \Vert \tilde{\omega} \Vert_2 + \Vert \tilde{F} \Vert_{L^2 \cdot W^{1,2}}  \right) $$
and 
$$ \nabla \Delta A + \Delta A V -\nabla A w + A (\nabla \omega +F)= \mathrm{curl} B.$$
Then we extend $u$ by $\tilde{u}\in W^{2,2}(B_1)$ such that 
$$ \Vert \nabla^2 \tilde{u} \Vert_{L^2 (B_1)} + \Vert \nabla \tilde{u} \Vert_{L^4 (B_1)}  \leq 2 \left( \Vert \nabla^2 u\Vert_{L^2 (B_1\setminus B_r)} + \Vert \nabla u\Vert_{L^4 (B_1\setminus B_r)} \right).$$
We easily see that $\tilde{u}$ satisfies
$$\Delta(A\Delta \tilde{u})= div(K) +A f \hbox{ on }   B_1\setminus B_r , $$
with $K=2\nabla A \Delta \tilde{u}-\Delta A \nabla \tilde{u} + A w \nabla \tilde{u} + \nabla A(V\nabla \tilde{u})-A \nabla (V \nabla \tilde{u}) -B\nabla \tilde{u}\in L^{\frac{4}{3},1} (B_1)$ such that 
$$  \Vert K\Vert_{L^\frac{4}{3}}  \leq C\left(1 + \Vert \nabla^2 u\Vert_{L^2 (B_1\setminus B_r)} + \Vert \nabla u\Vert_{L^4 (B_1\setminus B_r)}   \right) $$
Then, we extend  $Af$ by $\tilde{f}\in L^p(B_1) $ such that 
$$\Vert \tilde{f} \Vert_{p} \leq 2\Vert Af\Vert_p .$$ Then let $ D\in W^{1,\frac{4}{3}}_0(B_1)$ which satisfies
$$\Delta D =  div(K) + \tilde{f} \hbox{ on } B_1.$$
Hence, thanks to the standard $L^p$-theory, there exists $C$ a positive constant independent of $r$, such that 
$$ \Vert D \Vert_{2,1} \leq C\left(\Vert K\Vert_{L^{\frac{4}{3},1}} + \Vert f\Vert_p\right).$$ 
Finally, thanks to lemma \ref{l2}, there exists $a,b\in \R^k$ and $C$ a positive constant independent of $r$, such that 
\beq
\label{pp}
\begin{split}
\left\Vert D-A\Delta u - a -\frac{b}{\vert x\vert^2} \right\Vert_{L^{2,1}\left(B_\frac{1}{2}\setminus B_{2r}\right)} &\leq C \left\Vert D-A\Delta u\right\Vert_2\\
& \leq C\left(1 + \Vert \nabla^2 \tilde{u}\Vert_{2} + \Vert K\Vert_{L^{\frac{4}{3},1}} + \Vert f\Vert_p   \right) .
\end{split} 
\eeq 
Hence we have 
$$div(A \nabla \tilde{u} ) =  a+ \frac{b}{\vert x\vert^2} + F \hbox { on } B_1\setminus B_r$$
with 
$$\Vert F \Vert_{L^{2,1}\left(B_\frac{1}{2}\setminus B_{2r}\right)} \leq  C\left( 1+ \Vert \nabla^2 u\Vert_{L^2 (B_1\setminus B_r)} + \Vert \nabla u\Vert_{L^4 (B_1\setminus B_r)} +\Vert f \Vert_p \right) .$$
Let us proceed to the following Hodge decomposition,
\beq
\label{se1}
A d \tilde{u} = d \alpha + d^* \beta,
\eeq
where $\alpha \in W^{1,2}_0(B_\frac{1}{2})$ and $\beta \in W^{1,2}(B_\frac{1}{2})$  satisfy
$$
\Delta \alpha = a+ \frac{b}{\vert x\vert^2} + F \hbox{ on } B_\frac{1}{2} \setminus B_{2r}
$$
and 
$$
\Delta \beta = d A\wedge d\tilde{u}- d\tilde{u} \wedge d A    \hbox{ on } B_\frac{1}{2}.
$$
From the one hand, we extend $F$ by $\tilde{F}\in W^{1,2}\left(B_\frac{1}{2}\right)$ such that 
$$\Vert \tilde{F} \Vert_{L^{2,1}\left(B_\frac{1}{2}\right)} \leq 2 \Vert F\Vert_{L^{2,1}}.$$
Then, let $\tilde{\alpha}\in W^{1,2}_0(B_\frac{1}{2})$ which satisfies 
$$\Delta \tilde{\alpha} =  \tilde{F} \hbox{ on } B_\frac{1}{2}.$$
Hence, thanks to the standard $L^p$-theory, there exists $C$ a positive constant independent of $r$, such that 
$$\Vert \nabla^2 \tilde{\alpha}\Vert_{2,1} \leq C \Vert F\Vert_{2,1}.$$ 
then, thanks to lemma \ref{l2}, there exists $C$ a positive constant independent of $r$, such that 
\beq
\label{se2}
\begin{split}
\left\Vert \nabla^T\nabla \left( \alpha -\tilde{\alpha} \right)\right\Vert_{L^{2,1}\left(B_\frac{1}{4} \setminus B_{4r}\right)} &\leq C  \Vert \nabla^2 (\alpha-\tilde{\alpha})\Vert_{2}\\
& \leq C  \left(\Vert F\Vert_{2,1} + \Vert \nabla^2 \beta \Vert_2 + \Vert \nabla A\nabla \tilde{u}\Vert_2 +  \Vert A\nabla \tilde{u}\Vert_2\right) .
\end{split}
\eeq
From the other hand, thanks to the standard-$L^p$-theory and Sobolev embeddings, we get
\beq
\label{se3}
\left\Vert \nabla^2  \beta  \right\Vert_{L^{2,1}\left(B_\frac{1}{4}\right)}  \leq C\left( 1+ \Vert \nabla^2 u\Vert_{L^2 (B_1\setminus B_r)} + \Vert \nabla u\Vert_{L^4 (B_1\setminus B_r)} \right) .
\eeq
Here we use the injection of $W^{1,2}$ into $L^{4,2}$.
 Finally, thanks to (\ref{se1}), (\ref{se2}), (\ref{se3}) and the fact that
$$\left\Vert \nabla^T \nabla u \right\Vert_{L^{2,1}} \leq C \left ( \left\Vert \nabla^T (A \nabla u) \right\Vert_{L^{2,1}} + \left\Vert \nabla^T A \nabla u) \right\Vert_{L^{2,1} }\right), $$
we get the desired estimate.\hfill$\square$

\section{Proof of theorem \ref{th1}} 
\label{s4}   
First we are going to separate $B_1$ in three parts: one where $u_n$ converges to a limiting solution, an other composed of some small  neighborhoods where the energy concentrates and where some bubbles blow and a third part which consists of some neck regions which join the first  two parts. This ''bubble-tree'' decomposition is by now classical, see \cite{Pa} for instance, hence we just sketch briefly how to proceed.\\ 
 
{\bf Step 1 : Finding the points of concentration}\\

Let  $\eps_0$ be  such that the $V,w,\omega$ and $F$ given by the section \ref{s1}  satisfy, thanks to (\ref{cc}), the hypothesis of theorem \ref{ereg}
as soon as $\Vert \nabla^2 u\Vert_2^2 + \Vert \nabla u\Vert_4^4 \leq \eps_0$. Then, thanks to (\ref{be}), we easily proved that there exist finitely many points $a^1,\dots, a^n$ where 
\beq
\label{cp}
\int_{B(a_i, r)}  \left(\vert \nabla^2 u\vert^2 + \vert \nabla u\vert^4  \right)\, dx \geq \eps_0 \hbox{ for all } r>0.
\eeq
Moreover, using theorem \ref{ereg}, we prove that  there exist  $f\in L^p(B_1,\R^k)$ and a $f$-approximate biharmonic maps $u_\infty\in W^{2, 2}(B_1,N)$ , such that, up to a subsequence,  
$$f_n \rightharpoonup f \hbox{ in }  L^p(B_1,\R^k) $$
and
$$\nabla u_n\rightarrow \nabla u_\infty \hbox{ in } W^{1,\bar{p}}_{loc}(B_1\setminus \{ a^1,\dots,a^n\}).$$ 

{\bf Step 2 : Blow-up around $a^i$ }\\

We choose $r_i>0$ such that 
$$\int_{B(a_i, r_i)}  \left(\vert \nabla^2 u_\infty\vert^2 + \vert \nabla u_\infty\vert^4  \right)\, dx \leq \frac{\eps_0}{4} .$$  
Then, we define a center of mass of  $B(a^i, r^i)$ with respect to $u_n$ in the following way
$$a_n^i= \left(\frac{\displaystyle\int_{B(a^i,r^i)} x^\alpha \vert \nabla^2 u_n\vert^2\, dx}{\displaystyle\int_{B(a^i,r^i)}\vert \nabla^é u_n\vert^2\, dx} \right)_{\alpha=1,\dots,4}.$$
Let $\lambda^i_n$ be a positive real such that 
 $$\int_{B(a^i_n, r^i)\setminus B(a^i_n, \lambda_n^i)} \left(\vert \nabla^2 u_n \vert^2 +\vert \nabla u_n \vert^4\right) \, dx =  \frac{\eps_0}{2}.$$  
Then we set $\widetilde{u}_n^i(x) = u_n( a^i_n+ \lambda_n^i x)$ and $N^i_n =B(a_n^i,r^i)  \setminus B(a_n^i, \lambda_n^i)$. Thanks to the conformal invariance, we easily see that
\[
\int_{B\left(0,\frac{r^i}{\lambda_n^i}\right)} \left(\vert \nabla^2 \tilde{u}_n^i\vert^2 +\vert \nabla \tilde{u}_n^i\vert^4\right) \ dx = \int_{B(a_n^i,r^i)} \left(\vert \nabla^2 u_n\vert^2 +\vert \nabla u_n\vert^4\right) \ dx\le M
\]
and $\tilde{u}_n^i$ still satisfies the equation of approximate biharmonic maps with the approximation $(\lambda^i_n)^4 \tilde{f}_n$  which goes to zero in $L^p$-norm.  Let  $a_i^j$ be  the possible points of concentration of $\widetilde{u}_n^i$ where
\beq
\label{cp2}
\int_{B(a_i^j, r)}  \left( \vert \nabla^2 \tilde{u}_n^i\vert^2+ \vert \nabla \tilde{u}_n^i\vert^4 \right)  \, dz \geq \eps_0 \hbox{ for all } r>0,
\eeq
Then, up of a subsequence, for each $i$,
$$\nabla \tilde{u}_k^i \rightarrow \nabla u_\infty^i  \hbox{ in } W^{1,\bar{p}}_{loc}(B_1\setminus \{ a^1_i,\dots,a^{n_i}_i\}) ,$$ 
where $u^i_\infty\in W^{2,2}(\R^4,N)$ is a biharmonic map.\\

{\bf Step 3 : Iteration}\\
 
 Two cases have to be considered separately:\\

Either $\tilde{u}_n^i$ is subject to some concentration phenomenon as (\ref{cp}), and then we find some new points of concentration, in such a case we apply step $2$ to our new concentration points. Or, $\widetilde{u}_n^i$  converges in  $W^{2,\bar{p}}_{loc}(\R^4)$ to a non trivial biharmonic map.\\

Of course this process has to stop, since we are assuming a uniform bound on $\Vert \nabla^2 u_n\Vert_2 + \Vert \nabla^2 u_n\Vert_4 $  and  each step is consuming at least the energy of a non trivial biharmonic map which is bounded from below thanks to the energy gap proved in theorem \ref{ereg}.\\ 

{\bf Analysis of a neck region: }\\

A neck region is an annullar region which is a  union of a finite number of annuli $N_n^i= B\left(a_n^i,\mu^i_n \right )\setminus B\left (a_n^i, \lambda_n^i\right)$ such that 
$$\lim_{k\rightarrow +\infty} \mu^i_n=0,$$
$$\lim_{k\rightarrow +\infty} \frac{\lambda_n^i} {\mu^i_n}=0,$$
and 
\beq
\label{aa1}
\int_{N_n^i} \left(\vert \nabla^2 u_n \vert^2 +\vert \nabla u_n \vert^4\right) \, dx \leq \frac{\eps_0}{2} 
\eeq
In order to prove  theorem~\ref{th1}, we start by proving  a weak estimate on the energy of gradient and the hessian in the region $N_n^i$.\\

First we remark that, for all $\eps>0$, there exists $r>0$ such that for all $\rho>0$ such that 
$$B_{2\rho}(a_n^i)\setminus B_\rho(a_n^i) \subset N_n^i(r)$$
where $N_n^i(r)= B\left(a_n^i,r\mu^i_n \right )\setminus B\left (a_n^i, \frac{\lambda_n^i}{r}\right)$, we have 
\beq \label{we1}
\int_{B_{2\rho}(a_n^i)\setminus B_\rho(a_n^i) } \left(\vert \nabla^2 u_n \vert^2 +\vert \nabla u_n \vert^4\right) \, dx \leq \eps\quad .
\eeq
If this is not  the case there would exist a sequence $\rho_n^i\rightarrow 0$ such that, up to a subsequence, $\hat{u}_n = u_n(a_n^i+ \rho_n^i z)$ convergesin$W^{2,\bar{p}}_{loc}(\R^4\setminus \{0\})$ to  $\hat{u}_\infty$, a non-trivial biharmonic map.
Using the fact that the $W^{2,2}$-norm of $\hat{u}_\infty$ is bounded, we deduce using Schwartz lemma that it has to be in fact a solution on the whole space.
Using  the energy gap proved in theorem \ref{ereg} we deduce that $\hat{u}_\infty$  is such that
\beq
\int_{N_k^i} \left(\vert \nabla^2 u_\infty \vert^2 +\vert \nabla u_\infty \vert^4\right) \, dx \geq \eps_0 ,
\eeq
which contradicts (\ref{aa1}).\\

Then for all $\eps>0$, there exists $r>0$ such that 
\beq
\label{Linfiny} 
\Vert \nabla^2 u_n \Vert_{L^{2,\infty}(N_n^i(r))} +  \Vert \nabla u_n \Vert_{L^{4,\infty}(N_n^i(r))} \leq \eps.
\eeq
Indeed, let $0<\eps<\eps_0$ and  $r>0$ such that, for all  
$$B_{2\rho}(a_n^i)\setminus B_\rho(a_n^i) \subset N_n^i(r)$$
we have 
\beq 
\int_{B_{2\rho}(a_n^i)\setminus B_\rho(a_n^i) } \left(\vert \nabla^2 u_n \vert^2 +\vert \nabla u_n \vert^4 \right) \, dx \leq \eps \quad.
\eeq 
Then,thanks to $\eps$-regularity in theorem \ref{ereg}, there exist $q>2$ and  $C$ a positive constant, independent of $r$ and $u$, such that for all $\rho>0$ such that 
$$ B_{2\rho}(a_n^i)\setminus B_\rho(a_n^i) \subset N_n^i\left(\frac{r}{2}\right),$$
and  $n$ big enough, we have 
\beq
\label{linf}
\begin{split}
\rho^{2-\frac{4}{q}} \Vert \nabla^2  u \Vert_{L^{q}(B_{2\rho}\setminus B_\rho)} + \rho^{1-\frac{2}{q}} \Vert \nabla u \Vert_{L^{2q}(B_{2\rho}\setminus B_\rho)}  &\leq C\left( \sqrt{\eps}+ (r\mu_i^n)^\frac{4(p-1)}{p}\vert f_n\vert^p \right)\\
&\leq C \sqrt{\eps} . 
\end{split}
\eeq
Let $\lambda>0$, $f(x)=\vert \nabla^2 u(x)\vert$ if $x\in N_n^i(\frac{r}{2})$ and $f=0$ otherwise. For any $\rho>0$, we  denote 
$$U(\lambda,\rho ) \equiv \{x\in B_{2\rho} \setminus B_\rho \hbox{ s.t. } f(x) >\lambda\} .$$
Thanks to (\ref{linf}), we have  
$$\lambda^q \vert U(\lambda,\rho)\vert \leq C^r\eps^\frac{q}{2} \rho^{4-2q}.$$
Let $k\in\Z$ and $j\geq k$, we apply the previous inequality with $\rho= 2^{-j}\lambda^{-1} $ and we sum for $j\geq k$, which gives
$$\lambda^2 \vert \{x\in R^4  \setminus B_{2^k\lambda^{-1 }}  \hbox{ s.t. } f(x) >\lambda\}\vert \leq C 2^{-k(4-2q)}\eps^\frac{r}{2} \rho^{4-2q}.$$
Hence, for any $k\in \Z$, we have 
$$\lambda^2 \vert \{x\in R^4   \hbox{ s.t. } f(x) >\lambda\}\vert \leq C\left( 2^{-k(4-2q)} \eps^\frac{q}{2} +  2^{4k}\right).$$
Taking $2^{4k} \sim\eps^\frac{q}{2} $ we have
$$\Vert \nabla^2 u_n \Vert_{L^{2,\infty}(N_n^i(r))}  \leq C \eps^\frac{q}{4},$$ 
We prove a similar inequality for $\Vert \nabla u_n \Vert_{L^{4,\infty}}$, and then we have  (\ref{Linfiny}).\\

Finally using theorem~\ref{EN} and the duality for Lorentz spaces, we see that, for all $\eps>0$, there exists $r>0$ such that 
\beq
\label{L2L4} 
\Vert \nabla^T(\nabla u) \Vert_{L^{2}(N_k^i(r))}  \leq \eps 
\eeq
Then using the Poho\v{z}aev identity (\ref{pe}) for extrinsic biharmonic maps (reps. (\ref{pi}) for intrinsic biharmonic maps) and the fact that the convergence is strong on the boundary of a neck region, we get  that for all $\eps>0$, there exists $r>0$ such that 
\beq
\Vert \nabla^2 u \Vert_{L^{2}(N_k^i(r))} +  \Vert \nabla u \Vert_{L^{4}(N_k^i(r))} \leq \eps .
\eeq
Which achieves the proof of  theorem~\ref{th1}.\hfill$\square$\\

Following step by step the proof of theorem \ref{th1}, we can prove the following theorem about the angular energy quantization of solution of fourth order elliptic system in the form of Lamm Riviere, \cite{LR}.

\begin{thm}
\label{th2}
Let $V_n \in W^{1,2}(B_1, \mathcal{M}_k \otimes \Lambda^{1} \R^4)$, $w_n \in L^2(B_1, \mathcal{M}_k)$, $\omega_n \in L^2(B_1,\mathrm{s}o_k)$, $F_n \in L^2\cdot W^{1,2}(B_1,  \mathcal{M}_k \otimes \Lambda^{1} \R^4)$ and $u_n\in W^{2,1}(B_1,\R^n)$ be a sequence of solutions of
\beq
\Delta^2 u_n =   \Delta (V_n \nabla u_n) + div( w_n \nabla u_n) + \nabla \omega_n \nabla u_n + F_n \nabla u_n  ,
\eeq 
with bounded energy, i.e.
\beq
\Vert \nabla^2 u_n\Vert_2 +\Vert \nabla u_n\Vert_4 + \Vert V_n \Vert_{W^{1,2}} + \Vert w_n \Vert_2 +  \Vert \omega_n \Vert_2+  \Vert F_n \Vert_{L^2\cdot W^{1,2}}   \leq M.
\eeq
Then there exists $V_\infty \in W^{1,2}(B_1, \mathcal{M}_k \otimes \Lambda^{1} \R^4)$, $w_\infty \in L^2(B_1, \mathcal{M}_k)$, $\omega_\infty \in L^2(B_1,\mathrm{s}o_k)$, $F_\infty \in L^2\cdot W^{1,2}(B_1,  \mathcal{M}_k \otimes \Lambda^{1} \R^4)$ and let $u_\infty\in W^{2,1}(B_1,\R^n)$ a solution of 
$$\Delta^2 u_\infty= \Delta (V_\infty \nabla u_\infty) + div( w_\infty \nabla u_\infty) + \nabla \omega_\infty \nabla u_\infty + F_\infty \nabla u_\infty \hbox{ on }B_1,$$
$l\in\N^*$ and
\begin{enumerate}
\item $\theta^1,\dots, \theta^l$ a family of solutions to system of the form
$$\Delta^2 \theta^i = \Delta (V_\infty^i \nabla  \theta^i ) + div( w_\infty^i  \theta^i ) + \nabla \omega_\infty^i \nabla  \theta^i  + F_\infty^i \nabla  \theta^i  \mbox{ on } \R^4
$$
where $V_\infty^i \in W^{1,2}(\R^4, \mathcal{M}_k \otimes \Lambda^{1} \R^4)$, $w_\infty^i \in L^2(\R^4, \mathcal{M}_k)$, $\omega_\infty^i \in L^2(\R^4,\mathrm{s}o_k)$ and\\ $F_\infty^i \in L^2\cdot W^{1,2}(\R^4,  \mathcal{M}_k \otimes \Lambda^{1} \R^4)$,
\item $a_n^1,\dots,a_n^l$ a family of converging sequences of points of $B_1$,
\item  $\lambda_n^1,\dots,\lambda_n^l$ a family of  sequences of positive reals converging all to zero,
\end{enumerate}
such that, up to a subsequence,  
\be
\begin{split}
&V_n \rightharpoonup V_\infty \hbox{ in }  W^{1,2}_{loc}(B_1, \mathcal{M}_k \otimes \Lambda^{1} \R^4),\\
&w_n \rightharpoonup w_\infty \hbox{ in }  L^2_{loc}(B_1, \mathcal{M}_k),\\
&\omega_n \rightharpoonup \omega_\infty \hbox{ in }  L^2_{loc} (B_1,\mathrm{s}o_k),\\
&F_n \rightharpoonup F_\infty \hbox{ in }L^2_{loc}\cdot W^{1,2}_{loc}(B_1,  \mathcal{M}_k \otimes \Lambda^{1} \R^4),\\
& u_n\rightarrow u_\infty \hbox{ on } W^{2,2}_{loc}(B_1 \setminus \{ a^1_\infty,\dots,a^l_\infty\})  \\
&\hbox{ and }\\
& \left\Vert \left\langle\nabla \left( \nabla \left(u_n -u_\infty -\sum_{i=1}^l \theta^i_k \right)\right),X_n\right\rangle \right\Vert_{L^2_{loc}(B_1)} + \left\Vert \left\langle\nabla  \left(u_n -u_\infty -\sum_{i=1}^l \theta^i_k \right),X_n\right\rangle \right\Vert_{L^4_{loc}(B_1)} \rightarrow 0,
\end{split}
\ee
where $\omega^i_n=\omega^i(a_n^i+\lambda_n^i \, . \, )$ and  $X_n$ is any vector field whose image is in  $(\nabla d_n)^\bot $ with $\ds d_n= \min_{1\leq i\leq l } (\lambda_n^i + d(a_n^i,\, .\,))$.
\end{thm}

\section{A lemma about harmonic maps on an annular regions}
\label{s5}
\begin{lemma}
\label{l2}
Let $0<r<\frac{1}{8}$ and  $u\in W^{1,2}(B_1\setminus B_r)$ be a harmonic function such that
\be
\begin{split}
&\int_{\partial B_1} u \, d\sigma =0,\\
&\int_{\partial B_r} u \, d\sigma =0.
\end{split}
\ee
Then there exists $C$ a positive constant independent of $r$ and $u$ such that 
$$\Vert u\Vert_{L^{2,1}\left(B_\frac{1}{2}\setminus B_{2r}\right)} \leq C \Vert u \Vert_2 $$
and
$$\Vert\nabla^T  \nabla u\Vert_{L^{2,1}\left(B_\frac{1}{2}\setminus B_{2r}\right)} \leq C \Vert\nabla^T \nabla u \Vert_2 .$$
\end{lemma}

{\it Proof of lemma \ref{l2}:}\\

Since $u$ is harmonic, it can be decomposed with respect to the spherical harmonics as follows
\beq
\label{lap2}
u =\sum_{l=1}^{+\infty}\sum_{k=1}^{N_l} \left(d^l_k r^l +d^{-l}_kr^{-l-2}\right)\phi_k^l ,
\eeq 
where $(\phi_k^l)_{l,k}$ are a $L^2$-basis of  eigenfunction of the Laplacian on $S^3$. In particular we get 
$$\Delta \phi_k^l = - l(l+2) \phi_k^l \hbox{  on  } S^3.$$
Thanks to this equation, $L^p$-theory for singular operators gives the existence of a positive constant  $C$, independent of $l$ such that 
$$\Vert \phi_k^l\Vert_{\infty}\leq C((l(l+2))^2.$$ 
Moreover we know that  $N_l$, the dimension of the eigenspace associated to $-l(l+2)$, is equal to $(l+1)^2$.
Hence, computing  the $L^{2}$-norm and $L^{2,1}$-norm of the function $f_j:x\mapsto \vert x \vert^{j}$, we get 
$$\Vert f_j \Vert_2 \geq \frac{r^{2+j}}{2\sqrt{-2j-4}} \hbox {if }j<-2$$
$$\Vert f_j \Vert_2 \geq \frac{1}{2\sqrt{2j+4}} \hbox {if }j\geq0$$
$$\Vert f_j \Vert_{L^{2,1}\left(B_\frac{1}{2}\setminus B_{2r}\right)} \leq  (2r)^{2+j} \hbox {if }j <-2$$
$$\Vert f_j \Vert_{L^{2,1}\left(B_\frac{1}{2}\setminus B_{2r}\right)} \leq  \left(\frac{1}{2}\right)^{\frac{3j}{4}+1 } \hbox {if }j\geq 0$$
where $C$ is independent of $j$.\\

Then
\be
\begin{split}
\Vert u\Vert_{L^{2,1}\left(B_\frac{1}{2}\setminus B_r\right)} &\leq C  \sum_{l=1}^{+\infty}\sum_{k=1}^{N_l} \left(d^l_k \left(\frac{1}{2}\right)^{\frac{3l}{4}+1 }  +d^{-l}_k  (2r)^{-l} \right)((l(l+2))^2\\
& \leq C\left( \left( \sum_{l=1}^{+\infty}\sum_{k=1}^{N_l} (d^l_k)^2 \frac{1}{4(2l+4)}  \right)^\frac{1}{2} \left( \sum_{l=1}^{+\infty}\sum_{k=1}^{N_l}  4(2l+4)((l(l+2))^4 \left(\frac{1}{2}\right)^{\frac{3l}{2}+2 }   \right)^\frac{1}{2}\right. \\ 
&+\left.\left( \sum_{l=1}^{+\infty}\sum_{k=1}^{N_l} (d^{-l}_k)^2 \frac{r^{-2l}}{8l}  \right)^\frac{1}{2} \left( \sum_{l=1}^{+\infty}\sum_{k=1}^{N_l}   8l((l(l+2))^4 \left(\frac{1}{4}\right)^l  \right)^\frac{1}{2}\right)
\end{split}
\ee

Thanks to the fact that $N_l$, the dimension of the eigenspace associated to the eigenvalue $-l(l+2)$ of the Laplacian, is equal to $(l+1)^2$, we get the first estimate. The second identity is obtained in the same way. \hfill$\square$

\section{Poho\v{z}aev identities}
\label{sp}
In this section, we prove a Poho\v{z}aev identity for extrinsic and intrinsic biharmonic maps in order to rely the radial derivatives to the angular ones. First we multiply our equation by $x^k \partial_k u$ and we integrate by part.
\be
\begin{split}
\int_{B(0,r)} \left(x^k\partial_k u\right) (\Delta^2 u) \, dx &= -\int_{B(0,r)} \langle \nabla u  , \nabla (\Delta u)\rangle \, dx -\int_{B(0,r)} (x^k \partial_k \partial^i u)  (\partial_i(\Delta u)) \, dx\\
& +  \int_{\partial B(0,r)} \left( x^k \partial_k u\right)   \partial_\nu (\Delta u) \, d\sigma\\
&= 2 \int_{B(0,r)} (\Delta u)^2 \, dx + \int_{B(0,r)} x^k \partial_k ( \Delta u)   (\Delta u) \, dx\\
& +  \int_{\partial B(0,r)} \left(\left( r \partial_\nu u\right) \partial_\nu (\Delta u) -  \left(\partial_\nu u\right) (\Delta u)  -  r \left(\partial_\nu^2  u\right) (\Delta u)\right)  \, d\sigma\\
&=  \int_{\partial B(0,r)} \frac{r}{2}  (\Delta u)^2  \, d\sigma\\
& +  \int_{\partial B(0,r)}\left(\left( r \partial_\nu u\right)   \partial_\nu (\Delta u) -  \left(\partial_\nu u\right) (\Delta u)  -  r \left(\partial_\nu^2  u\right) (\Delta u)\right)  \, d\sigma
\end{split}
\ee
Using the fact that for an extrinsic harmonic maps we have $\Delta^2 u \bot T_u N $ almost everywhere, we get for all $r$ that
 
 \beq
 \label{p1}
 \int_{\partial B(0,r)}  \left(\frac{1}{2}(\Delta u)^2-   \left(\partial_\nu^2  u\right) \Delta u  + \left(  \partial_\nu u\right)   \partial_\nu (\Delta u) -\frac{1}{r}  \left(\partial_\nu u\right) (\Delta u)  \right)  \, d\sigma=0
\eeq
But 
$$\Delta u = \partial_\nu^2 u + \frac{3}{r}\partial_\nu u + \frac{1}{r^2} \Delta_{S^3}u .$$
Hence 
$$(\Delta u)^2 = (\partial_\nu^2 u)^2 + \frac{9}{r^2} (\partial_\nu u)^2 + \frac{1}{r^4} (\Delta_{S^3}u)^2+  \frac{6}{r} (\partial_\nu u)(\partial_\nu^2 u) +  \frac{2}{r^2} (\Delta_{S^3}u)(\partial_\nu^2 u)+  \frac{6}{r^3}(\partial_\nu u)(\Delta_{S^3}u).$$
From the one hand, we have

$$\frac{1}{2}(\Delta u)^2-   \left(\partial_\nu^2  u\right) \Delta u= -\frac{1}{2} \left(\partial_\nu^2  u\right)^2 +\frac{9}{2r^2} \left(\partial_\nu  u\right)^2 + \frac{1}{2r^4}\left(\Delta_{S^3} u\right)^2 + \frac{3}{r^3} \left(\partial_\nu  u\right)\left(\Delta_{S^3} u\right)$$
which gives
\beq
\label{p2}
\int_{B_R \setminus B_r } \left( \frac{1}{2}(\Delta u)^2- \left(\partial_\nu^2  u\right) \Delta u \right)\, dx = \int_{B_R \setminus B_r } \left( -\frac{1}{2} \left(\partial_\nu^2  u\right)^2 +\frac{9}{2r^2} \left(\partial_\nu  u\right)^2 + \frac{1}{2r^4}\left(\Delta_{S^3} u\right)^2 + \frac{3}{r^3} \left(\partial_\nu  u\right)\left(\Delta_{S^3} u\right)\right) \, dx
\eeq
From the other hand 
$$ \left(  \partial_\nu u\right)   \partial_\nu (\Delta u) -\frac{1}{r}  \left(\partial_\nu u\right) (\Delta u) =   \left(  \partial_\nu u\right)   \left(  \partial_\nu^3 u\right) + \frac{2}{r}  \left(  \partial_\nu u\right)  \left(  \partial_\nu^2 u\right) - \frac{6}{r}  \left(  \partial_\nu u\right) ^2 +\frac{1}{r^2}  \left(  \partial_\nu \Delta_{S^3 }u\right)  \left(  \partial_\nu u\right) -\frac{3}{r^3}  \left(\Delta_{S^3 }u\right)  \left(  \partial_\nu u\right)$$ 

Integrating by part, we get 
\beq
\label{p3}
\begin{split}
\int_{B_R \setminus B_r }  \left( \left(  \partial_\nu u\right)   \partial_\nu (\Delta u) -\frac{1}{r}  \left(\partial_\nu u\right) (\Delta u)\right)  \, dx &= \int_{B_R \setminus B_r }   \left( \left(  \partial_\nu u\right)   \left(  \partial_\nu^3 u\right) + \frac{2}{r}  \left(  \partial_\nu u\right)  \left(  \partial_\nu^2 u\right) - \frac{6}{r}  \left(  \partial_\nu u\right) ^2\right)\, dx \\
&+ \int_{B_R \setminus B_r }   \left( \frac{1}{r^2}  \left(  \partial_\nu \Delta_{S^3 }u\right)  \left(  \partial_\nu u\right) -\frac{3}{r^3}  \left(\Delta_{S^3 }u\right)  \left(  \partial_\nu u\right) \right)\, dx \\
&= \int_{\partial(B_R \setminus B_r )}  \left(  \partial_\nu u\right)   \left(  \partial_\nu^2 u\right)  \, d\sigma \\
&+ \int_{B_R \setminus B_r } \left( -\frac{1}{2r} \left(  \partial_\nu   \left(  \partial_\nu u\right)^2 \right)- \left(  \partial_\nu^2 u\right)^2   - \frac{6}{r}  \left(  \partial_\nu u\right) \right)\, dx\\
&+ \int_{B_R \setminus B_r }   \left( \frac{1}{r^2}  \left(  \partial_\nu \Delta_{S^3 }u\right)  \left(  \partial_\nu u\right) -\frac{3}{r^3}  \left(\Delta_{S^3 }u\right)  \left(  \partial_\nu u\right) \right)\, dx \\
&= \int_{\partial(B_R \setminus B_r )}  \left(\left(  \partial_\nu u\right)   \left(  \partial_\nu^2 u\right)  - \frac{1}{2r}  \left(  \partial_\nu u\right)^2\right) \, d\sigma \\
&-\int_{B_R \setminus B_r }  \left( \left(  \partial_\nu^2 u\right)^2   + \frac{5}{r^2}  \left(  \partial_\nu u\right)^2\right) \, dx\\
&+ \int_{B_R \setminus B_r }   \left( \frac{1}{r^2}  \left(  \partial_\nu \Delta_{S^3 }u\right)  \left(  \partial_\nu u\right) -\frac{3}{r^3}  \left(\Delta_{S^3 }u\right)  \left(  \partial_\nu u\right) \right)\, dx 
\end{split}
\eeq
Finally, thanks to (\ref{p1}), (\ref{p2}) and (\ref{p3}), we have 

\beq
\label{pe}
\begin{split}
 \int_{B_R \setminus B_r } \left( \frac{3}{2} \left(  \partial_\nu^2 u\right)^2   + \frac{1}{2r^2}  \left(  \partial_\nu u\right)^2\right) \, dx &=  \int_{B_R \setminus B_r } \left(  \frac{1}{2r^4}\left(\Delta_{S^3} u\right)^2\right) \, dx\\
&+ \int_{B_R \setminus B_r }   \left( \frac{1}{r^2}  \left(  \partial_\nu \Delta_{S^3 }u\right)  \left(  \partial_\nu u\right)  \right)\, dx \\
&+\int_{\partial(B_R \setminus B_r )}  \left(\left(  \partial_\nu u\right)   \left(  \partial_\nu^2 u\right)  - \frac{1}{2r}  \left(  \partial_\nu u\right)^2\right) \, d\sigma 
\end{split}
\eeq
Since the equations of extrinsic and intrinsic biharmonic maps differ only by  $P(u)\left(B(u) (\nabla u, \nabla u )  \nabla_u B(u) (\nabla u, \nabla u ) \right)+ 2 B(u)  (\nabla u , \nabla u) B(u) (\nabla u, \nabla P(u))$, we multiply this term by $x^k \partial_k u$ which gives
\be
\begin{split}
 & x^k\partial_k u \left(P(u)\left(B(u) (\nabla u, \nabla u )  \nabla_u B(u) (\nabla u, \nabla u ) \right)+2 B(u)  (\nabla u , \nabla u) B(u) (\nabla u, \nabla P(u))\right) \\
 &= B(u) (\nabla u, \nabla u )  \nabla_{x^k\partial_k u} B(u) (\nabla u, \nabla u ) +2 B(u)  (\nabla u , \nabla u) B(u) (\nabla u, \nabla (x^k\partial_k u )\\
 &= x^k\partial_k \left( \frac{\left\vert B(u)(\nabla u, \nabla u  )\right\vert^2}{2}\right) + 2 \left\vert B(u)(\nabla u, \nabla u  )\right\vert^2\\
 &=  \frac{1}{\vert x\vert^3} \frac{\partial}{\partial\nu }\left[ \frac{r^4}{2} \left\vert B(u)(\nabla u, \nabla u  )\right\vert^2\right].
 \end{split}
\ee
Then integrating, we get the following Poho\v{z}dev identity for intrinsic biharmonic maps
\beq
\label{pi}
\begin{split}
 \int_{B_R \setminus B_r } \left( \frac{3}{2} \left(  \partial_\nu^2 u\right)^2   + \frac{1}{2r^2}  \left(  \partial_\nu u\right)^2\right) \, dx &=  \int_{B_R \setminus B_r } \left(  \frac{1}{2r^4}\left(\Delta_{S^3} u\right)^2\right) \, dx\\
&+ \int_{B_R \setminus B_r }   \left( \frac{1}{r^2}  \left(  \partial_\nu \Delta_{S^3 }u\right)  \left(  \partial_\nu u\right)  \right)\, dx \\
&+\int_{\partial(B_R \setminus B_r )}  \left(\left(  \partial_\nu u\right)   \left(  \partial_\nu^2 u\right) \right.\\
&\left. - \frac{1}{2r}  \left(  \partial_\nu u\right)^2  - \frac{r}{2}\left\vert B(u)(\nabla u, \nabla u  )\right\vert^2 \right) \, d\sigma 
\end{split}
\eeq

We also get a Poho\v{z}dev identity for the critical point of general functional, since
$$\int_{B_R\setminus B_r } (x^k \partial_k u)  H\left(\frac{\partial u}{\partial x_1},\frac{\partial u}{\partial x_2},\frac{\partial u}{\partial x_3},\frac{\partial u}{\partial x_4}\right) \, dx = \int_{B_R\setminus B_r } d\Omega\left(x^k \frac{\partial u}{\partial x_k} ,\frac{\partial u}{\partial x_1},\frac{\partial u}{\partial x_2},\frac{\partial u}{\partial x_3},\frac{\partial u}{\partial x_4}\right) \, dx=0 .$$

\end{document}